\documentclass{jctart08}
\usepackage{graphics, graphicx}     

\usepackage{cite}
\usepackage{amsmath}
\usepackage{amssymb}
\usepackage{amsthm}
\theoremstyle{plain}
\newtheorem{Theorem}{Theorem}[section] %
\newtheorem{Lemma}{Lemma}[section]

\theoremstyle{definition}
\newtheorem{Def}{Definition}[section]

\renewcommand{\thetable}{\thesection.\arabic{table}}

\addto\captionsrussian{\def\refname{References}}
\newenvironment{Proof} 
{\par\noindent{\it Proof of}} 
{\hfill$\vspace{5mm}\scriptstyle\blacksquare$} 

\numberwithin{equation}{section} 
\numberwithin{figure}{section} 
\numberwithin{table}{section} 

\begin{document}

\setcounter{page}{1}

\markboth{M.I. Isaev}{Exponential instability in the inverse scattering problem on the energy interval}

\title{Exponential instability in the inverse scattering problem on the energy interval}
\date{}
\author{ {\bf M.I. Isaev}\\
Moscow Institute of Physics and Technology,
141700 Dolgoprudny, Russia\\
Centre de Math\'ematiques Appliqu\'ees, Ecole Polytechnique,
91128 Palaiseau, France\\
e-mail: \tt{isaev.m.i@gmail.com}}

\maketitle
{\bf Abstract}
\begin{abstract}
	We consider the inverse scattering problem on the energy interval in three dimensions. 
	We are focused on stability and instability questions for this problem. In particular, we prove an exponential instability estimate
	which shows optimality of the logarithmic stability result of [Stefanov, 1990] (up to the value of the exponent).
\end{abstract}

\section{Introdution}
We consider the Schr\"odinger equation
\begin{equation}\label{eq} 
	-\Delta\psi  + v(x)\psi = E\psi, \ \ \   x \in \mathbb{R}^3,
\end{equation}
where 
\begin{equation}\label{ubivanie}
	 \begin{array}{l}\displaystyle
	 v  \text{ is real-valued, } \ \  v\in L^{\infty}(\mathbb{R}^3),   \\ \displaystyle
	 v(x) = O(|x|^{-3-\epsilon}), \ \ \ \ |x|\rightarrow \infty, \ \ \text{ for some $\epsilon >0 $.}
	 \end{array}
\end{equation}

Under conditions (\ref{ubivanie}), for any
$k\in \mathbb{R}^3 \setminus 0$ equation (\ref{eq}) with $E = k^2$ has a unique continuous solution
$\psi^+(x,k)$ with asymptotics of the form
\begin{equation}\label{Psi_scattering}
	\begin{array}{c} \displaystyle
		\psi^+(x,k) = e^{ikx} - 2\pi^2 \frac{e^{i|k||x|}}{|x|} f\left(\frac{k}{|k|},\frac{x}{|x|},|k|\right)+ o\left(\frac{1}{|x|}\right) \\\\
		\displaystyle
		\text{as } \ |x| \rightarrow \infty \left( \text{uniformly in } \frac{x}{|x|}\right),
	\end{array}
\end{equation}	 
	where $f(k/|k|,\omega,|k|)$ with fixed $k$ is a continuous function of $\omega \in S^{2}.$

The function $f(\theta,\omega,s)$  arising in (\ref{Psi_scattering}) is refered to as the scattering amplitude for the potential $v$
for equation (\ref{eq}). (For more information on direct scattering for equation (\ref{eq}), under condition (\ref{ubivanie}), see, for
example, \cite{Faddeev1974} and \cite{Newton1989}.)

It is well known that for equation (\ref{eq}), under conditions (\ref{ubivanie}), the scattering amplitude $f$
in its high-energy limit uniquely determines $\hat{v}$ on $\mathbb{R}^3$ , where
\begin{equation}
	\hat{v}(p) = (2\pi)^{-3}\int_{\mathbb{R}^3} e^{ipx}v(x)dx, \ \ \ p\in \mathbb{R}^3, 
\end{equation}
via the Born formula. As a mathematical theorem this result goes back to \cite{Faddeev1956} (see, for example, Section 2.1
of \cite{Newton1989} and Theorem 1.1 of \cite{Novikov2001} for details).

We consider the following inverse problem for equation (\ref{eq}).

{\bf Problem 1.1} Given $f$ on the energy interval $I$, find $v$.



In \cite{Henkin1987} it was shown that for equation (\ref{eq}), under the conditions (\ref{ubivanie}), for any $E > 0$
and $\delta > 0$ the scattering amplitude $f(\theta,\omega,s)$ on $\{(\theta, \omega, s) \in S^2\times S^2 \times R_+,\   E \leq s^2 \leq E + \delta\}$ uniquely
determines $\hat{v}(p)$ on $\{p \in \mathbb{R}^3 \ | \ |p| \leq 2\sqrt{E} \}$.  This determination is based on solving linear integral
equations and on an analytic continuation. This result of \cite{Henkin1987} was improved in \cite{Novikov2001}.
 On the other hand, if $v$ satisfies (\ref{ubivanie}) and, in addition, is compactly
supported or exponentially decaying at infinity, then $\hat{v}(p)$ on $\{p \in \mathbb{R}^3 \ | \ |p| \leq 2\sqrt{E}\}$ uniquely determines
 $\hat{v}(p)$ on $\{p \in \mathbb{R}^3 \ | \ |p| > 2\sqrt{E}\}$ by an analytic continuation and, therefore,
 uniquely determines $v$ on $\mathbb{R}^3$.
 
 In the case of fixed energy and potential $v$, satisfying (\ref{ubivanie}) and, in addition, being compactly
supported or exponentially decaying at infinity, global uniqueness theorems and precise
reconstructions were given for the first time in \cite{Novikov1988}, \cite{Novikov1994}. 

An approximate but numerically efficient method for finding potential $v$ from the scattering amplitude $f$ in the case of fixed energy was devoloped in
\cite{Novikov2005}. Related numerical implementation was given in \cite{Alexeenko2008}. 

Global stability estimates for Problem 1.1 were given by Stefanov in \cite{Stefanov1990} (at fixed energy for compactly
supported potentials), 
see Theorem 2.1 in Section 2 of the present paper. 
In \cite{Stefanov1990}, using a special norm for
the scattering amplitude $f$, it was shown that the stability estimates for Problem 1.1 follow from the Alessandrini stability estimates of \cite{Alessandrini1988}  for the Gel'fand-Calderon inverse problem of finding potential $v$ in bounded domain from the Direchlet-to-Neumann map. The Alessandrini stability estimates  were recently improved by Novikov in \cite{Novikov2010}.

In the case of fixed energy, the Mandache results of \cite{Mandache2001} show that logarithmic 
stability estimates of Alessandrini of \cite{Alessandrini1988} and especially of Novikov of \cite{Novikov2010} are optimal (up to the value of the exponent). In \cite{IsaevDtN} studies of Mandache were extended to the 
case of Direchlet-to-Neumann map given on the energy intervals. Note also that
Mandache-type instability estimates for the elliptic inverse problem concerning
the determination of inclusions in a conductor by different kinds of boundary
measurements and the inverse obstacle acoustic scattering problems were given in \cite{Cristo2003}.

	 In the present work we apply to Problem 1.1 the approach of \cite{Mandache2001},\cite{IsaevDtN} and show that the Stefanov logarithmic 
stability estimates of \cite{Stefanov1990} are optimal (up to the value of the exponent). The Stefanov stability estimates and our instability result for Problem 1.1 are presented and discussed in Section 2. In Section 3 we prove some
 basic analytic properties of the scattering amplitude.
Finally, in Section 5 we prove the main result, using a ball packing and covering by ball arguments. 

\section{Stability and instability estimates}
In what follows we suppose 
\begin{equation}\label{supp_condition}
	\mbox{supp} \, v(x) \subset D = B(0,1), 
\end{equation}
where $B(x,r)$ is the open ball of radius $r$ centered at $x$. 
We consider the orthonormal basis of the spherical harmonics in $L^2(S^{2}) = L^2(\partial D)$:
\begin{equation}
		\{Y_{j}^{p} : j \geq 0;\  1 \leq p \leq 2j+1 \}. 
\end{equation}
The notation $(a_{j_1 p_1 j_2 p_2} )$ stands for a multiple sequence.
We will drop the subscript 
\begin{equation}
	0 \leq j_1,\  1 \leq p_1 \leq 2j_1+1 ,\  0 \leq j_2,\  1 \leq p_2 \leq 2j_2+1.  
\end{equation}
We expand function $f(\theta,\omega,s)$ in the basis $\{Y_{j_1}^{p_1}\times Y_{j_2}^{p_2}\}$:
\begin{equation}
	f(\theta,\omega,s) = \sum\limits_{j_1, p_1, j_2, p_2}a_{j_1 p_1 j_2 p_2}(s)Y_{j_1}^{p_1}(\theta)Y_{j_2}^{p_2}(\omega).
\end{equation} 
As in \cite{Stefanov1990} we use the norm
\begin{equation}\label{Norm_Stefanov}
	||f(\cdot,\cdot,s)||_{\sigma_1,\sigma_2} = 
	\left\{ \sum\limits_{j_1, p_1, j_2, p_2} 
	\left(\frac{2j_1+1}{es}\right)^{2j_1+2\sigma_1}
	\left(\frac{2j_2+1}{es}\right)^{2j_2+2\sigma_2}
	|a_{j_1 p_1 j_2 p_2}(s)|^2 
	\right\}^{1/2}.
\end{equation} 
If a function $f$ is the scattering amplitude for some potential $v\in L^\infty(D)$ supported in $B(0,\rho)$, where $0<\rho<1$, then
\begin{equation}\label{ineq_a}
	|a_{j_1 p_1 j_2 p_2}(s)| \leq C(s,||v||_{L^\infty(D)}) 	
	\left(\frac{es\rho}{2j_1+1}\right)^{j_1+3/2}
	\left(\frac{es\rho}{2j_2+1}\right)^{j_2+3/2}
\end{equation}
and,
therefore, $||f(\cdot,\cdot,s)||_{\sigma_1,\sigma_2}< \infty$, see estimates of Proposition 2.2 of \cite{Stefanov1990}.

\begin{Theorem}[see \cite{Stefanov1990}] 
	Let  $v_1, v_2$ be real-valued potentials such that 
	 $v_i \in L^{\infty}(D)\cap H^q(\mathbb{R}^3)$,  $\mbox{supp}\, v_i \subset B(0,\rho)$, 
	  $||v_i||_{L^{\infty}(D)} \leq N$ for  $i = 1,2$ and some $N>0$, $q>3/2$ and $0<\rho<1$.
		Let $f_1$ and $f_2$ denote the scattering amplitudes for $v_1$ and $v_2$, respectively, in
		the framework of equation (\ref{eq}) with $E = s^2,$ $s>0$, then
		\begin{equation}\label{estimation_stefanov}
			||v_1-v_2||_{L^{\infty}(D)} \leq c(N,\rho)\phi_\delta(||f_1(\cdot,\cdot,s)-f_2(\cdot,\cdot,s)||_{3/2,-1/2}),
		\end{equation}
	where  $\phi_\delta(t) = (-\ln t)^{-\delta}$ for some fixed $\delta$, where, in particular, $0<\delta<1$, and for sufficiently small $t>0$. 
	
\end{Theorem}
The main result of the present work is the following theorem.
\begin{Theorem}\label{main}
	  For the interval $I=[s_1,s_2]$, such that $s_1>0$, and for any
	$m>0$, $\alpha > {2m}$ and any real $\sigma_1,\sigma_2$ there are constants $\beta>0$ and $N>0$, such that for any 
	$v_0 \in C^m(D)$ with $||v_0||_{L^{\infty}(D)}\leq N$, $\mbox{supp}\, v_0 \subset B(0,1/2)$ and any $\epsilon \in (0,N)$, there are real-valued
	potentials $v_1,v_2 \in C^m(D)$, also supported in $B(0,1/2)$, such that
	\begin{equation}\label{maineq}
		\begin{array}{l}
			\displaystyle
				\sup_{s\in I} \left(||f_1(\cdot,\cdot,s)-f_2(\cdot,\cdot,s)||_{\sigma_1,\sigma_2} \right)
				\leq \exp\left( -\epsilon^{-\frac{1}{\alpha}}\right),\\
			\displaystyle
				||v_1-v_2||_{L^\infty(D)} \geq \epsilon, \\
			\displaystyle
				||v_i-v_0||_{L^\infty(D)} \leq \epsilon, \ \ \ \ i=1,2,\\
			\displaystyle
				||v_i-v_0||_{C^m(D)} \leq \beta, \ \ \ \ i=1,2,
		\end{array}
	\end{equation}
	where $f_1, f_2$ are the scattering amplitudes for $v_1, v_2$, respectively, for equation (\ref{eq}).
\end{Theorem}

\noindent
{\bf Remark 2.1.}  
In the case of fixed energy $s_1 = s_2$ we can replace the condition $\alpha > {2m}$ in Theorem \ref{main} by  $\alpha > {5m}/3$.

\noindent
{\bf Remark 2.2.}  We can allow $\beta$ to be arbitrarily small in Theorem \ref{main} if we 
 	require $\epsilon \leq \epsilon_0$ and replace the right-hand side in the first inequality in (\ref{maineq})
 	by $\exp(-c\epsilon^{-\frac{1}{\alpha}} )$, with $\epsilon_0 > 0$ 	and $c > 0$ depending on $\beta.$
 	
\noindent
{\bf Remark 2.3.}  	Note that Theorem \ref{main} and Remark 2.1 imply, in particular, that for any real $\sigma_1$ and $\sigma_2$ the estimate
 	\begin{equation}
 		||v_1-v_2||_{L^{\infty}(D)} \leq \tilde{c}(N,\rho,m, I)\sup_{s \in I}\phi_\delta(||f_1(\cdot,\cdot,s)-f_2(\cdot,\cdot,s)||_{\sigma_1,\sigma_2})
	\end{equation}
	 can not hold
 	with $\delta > 2m$ in the case of the scattering amplitude given on the energy interval and 
 	with $\delta > 5m/3$ in the case of fixed energy.
 	Thus Theorem \ref{main} and Remark 2.1 show optimality of
 	the Stefanov logarithmic stability result (up to the value of the exponent). 
 	
\noindent
{\bf Remark 2.4.} 	
 		A disadvantage of estimate (\ref{estimation_stefanov}) is that
\begin{equation}\label{alpha_A}
	\delta < 1 \text{ even if $m$ is very great.} 
\end{equation} 
Apparently, proceeding from results of \cite{Novikov2010}, it is not difficult to improve estimate (\ref{estimation_stefanov}) for 
\begin{equation}
	\delta = m + o(m) \text{ as } m \rightarrow \infty.
\end{equation}
\section{Some basic analytic properties of the scattering amplitude}
Consider the solution $\psi^+(x,k)$ of equation \ref{eq}, see formula (\ref{Psi_scattering}). We have that
\begin{equation}
	\psi^+(x,k) = e^{ikx}\mu^+(x,\theta,s),
\end{equation}
where $\theta \in S^2$, $k = s\theta$ and $\mu^+(x,\theta,s)$ solves the equation
\begin{equation}\label{eq_mu}
\mu^{+}(x,\theta,s) = 1 -  \int_{\mathbb{R}^3} G^+(x,y,s) e^{-is\theta(x-y)} v(y)
\mu^+(y,\theta,s)dy,
\end{equation}
where
\begin{equation}
G^+(x,y,s) = \frac{e^{is|x-y|}}{4\pi|x-y|}.
\end{equation}
We suppose that condition (\ref{supp_condition}) holds and, in addition, for some $h>0$ we have that
\begin{equation}
	|\mbox{Im}\,s| \leq h, 
\end{equation}
\begin{equation}\label{condition_normv}
	c_1(h,D)||v||_{L^{\infty}(D)} \leq 1/2, 
\end{equation}
where $D = B(0,1)$,
\begin{equation}
	c_1(h,D) = \sup_{x \in D} \int_D \frac{e^{2h|x-y|}}{4\pi|x-y|} dy.
\end{equation}
Then, in particular,
\begin{equation}
	\left|e^{-is\theta(x-y)}e^{is|x-y|}\right| \leq e^{2h|x-y|}.
\end{equation}
Solving (\ref{eq_mu}) by the method of succesive approximations in $L^{\infty}(D)$, we obtain that 
\begin{equation}\label{modul_mu}
	|\mu^+(x,\theta,s)| \leq \frac{1}{1-c_1||v||_{L^{\infty}(D)}}, \ \ \ \theta \in S^2, \ \   x \in D. 
\end{equation}
\begin{Lemma}\label{Lemma1}
	Let $a_{j_1 p_1 j_2 p_2}(s)$ denote coefficients $f(s,\theta,\omega)$ in the basis of the 
	spherical harmonics $\{Y_{j_1}^{p_1}\times Y_{j_2}^{p_2}\}$, where $f$ is the scattering amplitude
	 for potential $v\in L^{\infty}(D)$ such that conditions (\ref{supp_condition}) and (\ref{condition_normv}) hold  for some $h>0$,
	\begin{equation}
		f(\theta,\omega,s) = \sum\limits_{j_1, p_1, j_2, p_2}a_{j_1 p_1 j_2 p_2}(s)Y_{j_1}^{p_1}(\theta)Y_{j_2}^{p_2}(\omega).
	\end{equation}
	Then $a_{j_1 p_1 j_2 p_2}(s)$ is holomorphic function in $W_h = \{ s \ | \ |\mbox{Im}\,s| \leq h \}$ and 
	\begin{equation}
		|a_{j_1 p_1 j_2 p_2}(s)| \leq c_2(h,D)  
		s \in W_h.
	\end{equation}
\end{Lemma}
\begin{Proof}{ \it Lemma \ref{Lemma1}.}
	We start with the well-known formula
	\begin{equation}\label{formula_mu}
		f(\theta,\omega,s) = \frac{1}{(2\pi)^3}\int_{\mathbb{R}^3} e^{is(\theta-\omega)x}v(x)\mu^+(x,\theta,s)dx.
	\end{equation}
	Note that, since $\theta,\omega \in S^2$, 
	\begin{equation}
		|e^{is(\theta-\omega)x}| \leq e^{2|Im\,s||x|}. 
	\end{equation}
	Combining it with (\ref{supp_condition}), (\ref{condition_normv}), (\ref{modul_mu}) and  (\ref{formula_mu}) we obtain that 
	\begin{equation}
		|f(\theta,\omega,s)| \leq \tilde{c}_2(h,D) \text{ for } s \in W_h. 
	\end{equation}
	Using also that
	\begin{equation}
		a_{j_1 p_1 j_2 p_2}(s) = \int_{S^2 \times S^2}
		f(\theta,\omega,s)Y_{j_1}^{p_1}(\theta)Y_{j_2}^{p_2}(\omega) d\theta d\omega
	\end{equation}
	we obtain the result of Lemma \ref{Lemma1}.
\end{Proof}
\section{A fat metric space and a thin metric space}

\begin{Def}
	Let $(X,dist)$ be a metric space and $\epsilon > 0$. We say that a set $Y \subset X$ is an $\epsilon$-net
	for $X_1 \subset X$ if for any $x\in X_1$ there is $y \in Y$ such that $dist(x,y)\leq\epsilon.$ We call 
	$\epsilon$-entropy of the set $X_1$ the number $\mathcal H_\epsilon (X_1) := \log_2
	\min\{|Y|: Y$  is an $\epsilon$-net fot $X_1 \}.$ 
	
	A set $Z \subset X$ is called $\epsilon$-discrete if for any distinct $z_1, z_2 \in Z$, we have 		  
	$dist(z_1,z_2)\geq\epsilon$. We call $\epsilon$-capacity of the set $X_1$ the number $\mathcal C_\epsilon :=
	 \log_2	\max\{|Z|: Z \subset X_1 $ and $Z$ is $\epsilon$-discrete$\}.$
\end{Def}

	The use of $\epsilon$-entropy and $\epsilon$-capacity  to derive properties of mappings between metric
spaces goes back to Vitushkin and Kolmogorov (see \cite{Kolmogorov1959} and references therein). One notable
application was Hilbert’s 13th problem (about representing a function of several variables as
a composition of functions of a smaller number of variables). In essence, Lemma \ref{Lemma2} and Lemma \ref{Lemma3} are parts of the Theorem XIV and the Theorem XVII in \cite{Kolmogorov1959}. 

\begin{Lemma}\label{Lemma2}
	Let $d \geq 2$ и $m > 0$. For $\epsilon, \beta > 0$, consider the real metric space
	$$
		X_{m\epsilon\beta} = \{v \in C^{m}(\mathbb{R}^{d})\ |\   \mbox{supp}\, v \subset B(0,1/2),\  
		||v||_{L^\infty(\mathbb{R}^{d})} \leq \epsilon,\ 	||v||_{C^m(\mathbb{R}^{d})}\leq \beta\}
	$$
	with the metric induced by $L^{\infty}$. Then there is $\mu>0$ such that for any $\beta > 0$ and 			
	$\epsilon \in (0, \mu\beta)$, there is an
	$\epsilon$-discrete set  $Z \subset X_{m\epsilon\beta}$ with at least $\exp\Big(2^{-d-1}(\mu\beta/\epsilon)^{d/m}\Big)$ elements.
\end{Lemma}
\begin{Lemma}\label{Lemma3}
	For the interval $I = [a,b]$ and $\gamma>0$ consider the ellipse $W_{I,\gamma} \in \mathbb{C}$:
	\begin{equation}\label{rectangle}
		W_{I,\gamma} = \{ \frac{a+b}{2}+ \frac{a-b}{2}\cos z \ |\ |Im\, z| \leq \gamma \}.
	\end{equation}  	
	Then there is a constant $\nu = \nu(C,\gamma)>0$ such that for any $\delta\in(0,e^{-1})$ there is
	a $\delta$-net for the space of functions on $I$ with $L^{\infty}$-norm, having holomorphic continuation to $W_{I,\gamma}$ 
	 with module bounded above on $W_{I,\gamma}$ by the constant $C$, 
	with at most $\exp(\nu (\ln\delta^{-1})^2)$ elements.    
\end{Lemma}

\noindent
{\bf Remark 4.1.} In the case of $a=b$,  taking
\begin{equation}
	Y =\frac{\delta}{2}\mathbb{Z}\bigcap[-C,C] + i\cdot \frac{\delta}{2}\mathbb{Z}\bigcap[-C,C],
\end{equation}
we get $\delta$-net with at most $\exp(\nu \ln\delta^{-1})$ elements.

 Lemma \ref{Lemma2} and Lemma \ref{Lemma3} were also formulated and proved in \cite{Mandache2001} and  \cite{IsaevDtN}, respectively.

For the interval $I=[s_1,s_2]$ such that $s_1>0$ and real $\sigma_1, \sigma_2$ we introduce the Banach space
\begin{equation}
			X_{I, \sigma_1,\sigma_2} = \left\{ \Big(a_{j_1 p_1 j_2 p_2}(s)\Big)\ | \ 
			 \left\|\Big(a_{j_1 p_1 j_2 p_2}(s)\Big)\right\|_{X_{I, \sigma_1,\sigma_2}} < \infty \right\},
\end{equation}
where
\begin{equation}\label{NormXI}
	\left\|\Big(a_{j_1 p_1 j_2 p_2}(s)\Big)\right\|_{X_{I, \sigma_1,\sigma_2}} = 
\sup_{
\begin{array}{c}
\scriptscriptstyle
\vspace{-7.5mm}
s \in I
\\
\scriptscriptstyle
j_1, p_1, j_2, p_2
\end{array}
} 
\left(
\left(\frac{2j_1+1}{es}\right)^{j_1+\sigma_1}
	\left(\frac{2j_2+1}{es}\right)^{j_2+\sigma_2}
	|a_{j_1 p_1 j_2 p_2}(s)| \right).
\end{equation}   
We consider the scattering amplitude $f$ for some potential $v\in L^\infty(D)$ supported in $B(0,\rho)$, where $0<\rho<1$.
We identify in the sequel the scattering amplitude $f(s, \theta,\omega)$ with its matrix $\Big(a_{j_1 p_1 j_2 p_2}(s)\Big)$ 
in the basis of the 
	spherical harmonics $\{Y_{j_1}^{p_1}\times Y_{j_2}^{p_2}\}$ . We have that
\begin{equation}\label{ineq_f}
	\sup_{s\in I}||f(\cdot,\cdot,s)||_{\sigma_1,\sigma_2} \leq c_3
	\left\|\Big(a_{j_1 p_1 j_2 p_2}(s)\Big)\right\|_{X_{I, \tilde{\sigma}_1,\tilde{\sigma}_2}},
\end{equation}
where $\tilde{\sigma}_1 - \sigma_2 = \tilde{\sigma}_2 - \sigma_2 = 3$ and $c_3 = c_3(I) > 1$. 
We obtain (\ref{ineq_f}) from definitions (\ref{Norm_Stefanov}), (\ref{NormXI}) and by taking $c_3 > 1$ in a such a way that
\begin{equation}
	\sum\limits_{j_1, p_1, j_2, p_2}
	\left(\frac{2j_1+1}{es}\right)^{-3}
	\left(\frac{2j_2+1}{es}\right)^{-3}< c_3.
\end{equation}
For $h >0$ we denote by $\mathcal A_h$ the set of the matrices, corresponding to the scattering amplitudes for the potentials $v\in L^{\infty}(D)$ 
supported in $B(0,1/2)$ such that condition (\ref{condition_normv}) holds.
\begin{Lemma}\label{Lemma4}
For any $h>0$ and any real $\sigma_1, \sigma_2$, the set $\mathcal A_h$ belongs to $X_{I, \sigma_1,\sigma_2}$. In addition, there is a constant 
$\eta = \eta(I,h,\sigma_1,\sigma_2) > 0$ such that for any $\delta\in{(0,e^{-1})}$
 there is a $\delta$-net $Y$ for $\mathcal A_h$ in $X_{I, \sigma_1,\sigma_2}$ with at most 
 $\exp \left(\eta\left(\ln\delta^{-1}\right)^{6} \left(1 + \ln\ln\delta^{-1}\right)^2\right)$ elements.
\end{Lemma}
\begin{Proof} {\it Lemma \ref{Lemma4}}. 
	We can suppose that $\sigma_1, \sigma_2 \geq	0$ as the assertion is stronger in this case. If a function $f$ is the scattering 
	amplitude for some potential $v \in L^{\infty}(D)$ supported in $B(0, 1/2)$, we have from (\ref{ineq_a}) that
	\begin{equation}
		\left(\frac{2j_1+1}{es}\right)^{j_1+\sigma_1}
		\left(\frac{2j_2+1}{es}\right)^{j_2+\sigma_2}
		|a_{j_1 p_1 j_2 p_2}(s)|
		\leq
		c_4 \frac{(2j_1+1)^{\sigma_1}(2j_2+1)^{\sigma_2}}{ 2^{j_1+j_2}},
	\end{equation}
	where $c_4 = c_4(I,h) > 0.$ Hence, for any positive $\sigma_1$ and $\sigma_2$,
	\begin{equation}
		\left\|\Big(a_{j_1 p_1 j_2 p_2}(s)\Big)\right\|_{X_{I, \sigma_1,\sigma_2}} 
		\leq \sup_{j_1,j_2} \left( c_4 \frac{(2j_1+1)^{\sigma_1}(2j_2+1)^{\sigma_2}}{ 2^{j_1+j_2}}\right) < \infty
	\end{equation} 
	and so the first assertion of the Lemma \ref{Lemma4} is proved.
	
		Let $l_{\delta,\sigma_1,\sigma_2}$ be the smallest natural number such that $c_4(2l+1)^{\sigma_1+\sigma_2}2^{-l}<\delta$
		for any $l\geq l_{\delta,\sigma_1,\sigma_2}.$ Taking natural logarithm we have that 
		\begin{equation}
			- \ln c_4 - (\sigma_1+\sigma_2) \ln (2l+1) + l\ln2  > \ln{\delta^{-1}} \text{ for } l\geq l_{\delta,\sigma_1,\sigma_2}.
		\end{equation}
	 Using $\ln \delta^{-1} > 1$, we get that
		\begin{equation}\label{lsmax}
			l_{\delta,\sigma_1,\sigma_2} \leq C' \ln \delta^{-1},
		\end{equation}
		where the constant $C'$ depends only on  $h$, $\sigma_1$, $\sigma_2$ and $I=[s_1,s_2]$.
		 We take $W_I = W_{I,\gamma}$ of (\ref{rectangle}), where the constant $\gamma>0$ is such that  $W_I \subset \{ s \ | \ |\mbox{Im}\,s| \leq h \}$.
		If $\max(j_1,j_2)\leq l_{\delta, \sigma_1, \sigma_2}$, then we denote by $Y_{j_1 p_1 j_2 p_2}$ some $\delta_{j_1 p_1 j_2 p_2}$-net	
		from Lemma \ref{Lemma3} with the constant $C = c_2$, where the constant $c_2$ is from Lemma \ref{Lemma1} and 
		\begin{equation} 
			\delta_{j_1 p_1 j_2 p_2} = 	\left(\frac{es_1}{2j_1+1}\right)^{j_1+\sigma_1}
		\left(\frac{es_1}{2j_2+1}\right)^{j_2+\sigma_2} \delta. 
		\end{equation}
		Otherwise we take $Y_{j_1 p_1 j_2 p_2} = \{0\}.$
		We set
		\begin{equation}
			Y = \left\{
			\Big(a_{j_1 p_1 j_2 p_2}(s)\Big) \ |\ a_{j_1 p_1 j_2 p_2}(s)\in Y_{j_1 p_1 j_2 p_2} 
				 \right\}.
		\end{equation}
		For any $\Big(a_{j_1 p_1 j_2 p_2}(s)\Big) \in \mathcal A_h$ there is an element $\Big(b_{j_1 p_1 j_2 p_2}(s)\Big) \in Y$ such that
	
		\begin{equation}
			\begin{array}{c}\displaystyle
			\left(\frac{2j_1+1}{es}\right)^{j_1+\sigma_1}
		\left(\frac{2j_2+1}{es}\right)^{j_2+\sigma_2}|a_{j_1 p_1 j_2 p_2}(s)-b_{j_1 p_1 j_2 p_2}(s)| \leq\\\displaystyle
		\leq \left(\frac{2j_1+1}{es}\right)^{j_1+\sigma_1}
		\left(\frac{2j_2+1}{es}\right)^{j_2+\sigma_2}\delta_{j_1 p_1 j_2 p_2} \leq \delta
		\end{array}
		\end{equation}
		in the case of $\max(j_1,j_2)\leq l_{\delta, \sigma_1, \sigma_2}$ and
			\begin{equation}
			\begin{array}{c}\displaystyle
			\left(\frac{2j_1+1}{es}\right)^{j_1+\sigma_1}
		\left(\frac{2j_2+1}{es}\right)^{j_2+\sigma_2}|a_{j_1 p_1 j_2 p_2}(s)-b_{j_1 p_1 j_2 p_2}(s)| \leq\\\displaystyle
		\leq
		c_4 \frac{(2j_1+1)^{\sigma_1}(2j_2+1)^{\sigma_2}}{ 2^{j_1+j_2}} \leq 
				c_4\frac{\left(2\max(j_1,j_2)+1\right)^{\sigma_1+\sigma_2}}{2^{\max(j_1,j_2)}}
		< 
		\delta,
		\end{array}
		\end{equation}
		otherwise.
		
			It remains to count the elements of $Y$. We recall that $|Y_{j_1 p_1 j_2 p_2}| = 1$ in the case of $\max(j_1,j_2)> l_{\delta, \sigma_1, \sigma_2}$.
			 Using again the fact that $\ln\delta^{-1} \geq 1$ and (\ref{lsmax})
			we get in the case of $\max(j_1,j_2)\leq l_{\delta, \sigma_1, \sigma_2}$:  
			\begin{equation}\label{Yjpkq}
				|Y_{j_1 p_1 j_2 p_2}| \leq \exp(\nu (\ln\delta_{j_1 p_1 j_2 p_2}^{-1})^2) \leq 
				\exp \left(\nu' \left(\ln\delta^{-1}\right)^{2} \left(1 + \ln\ln\delta^{-1}\right)^2\right).
			\end{equation}
			We have that 
			$n_{\delta, \sigma_1, \sigma_2} \leq 
				l_{\delta, \sigma_1, \sigma_2}^2(2l_{\delta, \sigma_1, \sigma_2}+1)^{2}
				\leq (2l_{\delta, \sigma_1, \sigma_2}+1)^{4}
				,
			$ 
			where 
			$n_{\delta, \sigma_1, \sigma_2}$ is the number of four-tuples $(j_1, p_1, j_2, p_2)$ with $\max(j_1,j_2) \leq l_{\delta, \sigma_1, \sigma_2}$. 
			Taking $\eta$ to be big enough we get that
			\begin{equation}\label{allignY}
				\begin{aligned}
				|Y| &\leq \left(\exp \left(\nu' \left(\ln\delta^{-1}\right)^{2} \left(1 + \ln\ln\delta^{-1}\right)^2\right)\right)^{n_{\delta, \sigma_1, \sigma_2}} \\
						&\leq \exp\left(\nu' \left(\ln\delta^{-1}\right)^{2} \left(1 + \ln\ln\delta^{-1}\right)^2
						 {(1+2C'\ln\delta^{-1})^{4}}\right)\\
						&\leq \exp \left(\eta\left(\ln\delta^{-1}\right)^{6} \left(1 + \ln\ln\delta^{-1}\right)^2\right).
				\end{aligned}	
			\end{equation}
\end{Proof}

\noindent
{\bf Remark 4.2.} In the case of $s_1=s_2$,  taking into account  Remark 4.1 and using it in (\ref{Yjpkq}) and (\ref{allignY}),
we get $\delta$-net $Y$ with at most $\exp \left(\eta\left(\ln\delta^{-1}\right)^{5} \left(1 + \ln\ln\delta^{-1}\right)\right)$ elements.
\section{Proof of Theorem 2.2}

	We take $N$ such that condition (\ref{condition_normv}) holds for any $||v||_{L^{\infty}(D)}\leq 2N$ for some $h>0$. 
	By Lemma \ref{Lemma2}, the set $v_0 + X_{m\epsilon\beta}$ has an $\epsilon$-discrete subset $v_0 + Z$. 
	Since $\epsilon \in (0,N)$ we have that 
	the set $Y$
	constructed in Lemma \ref{Lemma4} is also $\delta$-net for the set of the matrices, corresponding to the scattering 
	amplitudes for the potentials  $v \in v_0+X_{m\epsilon\beta}$. 
	We take $\delta$ such that $2c_3\delta = \exp\left(-\epsilon^{-\frac{1}{\alpha}} \right)$, see (\ref{ineq_f}). 
	Note that inequalities of (\ref{maineq}) follow from
	\begin{equation}\label{ZplusY}
		|v_0 + Z|>|Y|,
	\end{equation}
	where the set $Y$ is constructed in Lemma \ref{Lemma4} with $\tilde{\sigma}_1 = \sigma_1+ 3$ and $\tilde{\sigma}_2 = \sigma_2+ 3$.
	In fact, if $|v_0 + Z|>|Y|$, then there are two potentials $v_1, v_2 \in v_0 + Z$ 
	with the matrices $\Big(a_{j_1 p_1 j_2 p_2}(s)\Big)$ and $\Big(b_{j_1 p_1 j_2 p_2}(s)\Big)$
	, corresponding to the scattering 
	amplitudes for them, 
	being in the same $X_{I, \sigma_1,\sigma_2}$-ball radius $\delta$ centered at a point of $Y$. Hence, using (\ref{ineq_f}) we get that 
	\begin{equation}
		\begin{aligned}
			\sup_{s\in I}||f_1(\cdot,\cdot,s) - f_2(\cdot,\cdot,s)||_{\sigma_1,\sigma_2} &\leq c_3
		\left\|\Big(a_{j_1 p_1 j_2 p_2}(s)\Big) - \Big(b_{j_1 p_1 j_2 p_2}(s)\Big)\right\|_{X_{I, \tilde{\sigma}_1,\tilde{\sigma}_2}}\leq\\
	 	&\leq2c_3\delta = \exp\left(-\epsilon^{-\frac{1}{\alpha}} \right).
	 	\end{aligned}
	\end{equation}
	It remains to find $\beta$ such that (\ref{ZplusY}) is fullfiled.	By Lemma \ref{Lemma4} for some 
	$\eta_\alpha = \eta_\alpha(I,\sigma_1,\sigma_2, \alpha) > 0$
	\begin{equation} \label{eqy}
		|Y|\leq 
		\exp\left(
		\eta\left( \ln(2c_3) +\epsilon^{-\frac{1}{\alpha}} \right)^{6} 
		\left(
			1 + \ln\left( \ln(2c_3) +\epsilon^{-\frac{1}{\alpha}} \right)
		\right)^2
		 \right)
		\leq 
		\exp\left(
			\eta_\alpha \epsilon^{-\frac{3}{m}}
		 \right)
		.
	\end{equation}
	Now we take 
	\begin{equation}\label{eqbeta}
		\beta > \mu^{-1} \max\left(N, \eta_\alpha^{m/3}2^{2m} \right). 
	\end{equation}
	This fulfils requirement $\epsilon<\mu\beta$ in Lemma \ref{Lemma2}, which gives
	\begin{equation}
		\begin{array}{c}\displaystyle
			|v_0+Z| = |Z| \geq \exp\Big(2^{-4}(\mu\beta/\epsilon)^{3/m}\Big)\stackrel{(\ref{eqbeta})}{>}\\\displaystyle
			>\exp\left(2^{-4}(\eta_\alpha^{m/3}2^{2m}/\epsilon)^{3/m}\right) \stackrel{(\ref{eqy})}{\geq} |Y|.
		\end{array}
	\end{equation}
This completes the proof of Theorem \ref{main}.

In the case of fixed energy $s_1 = s_2$, using Remark 4.2 in (\ref{eqy}), we can replace the condition $\alpha > 2m$ in Theorem \ref{main} by  $\alpha > 5m/3$.

\section*{Acknowledgments}
This work was fulfilled under the direction
of R.G.Novikov in the framework of an internship at Ecole Polytechnique.

\end{document}